%notes on lattice polytopes
\input amstex
%This will run under Ams-Tex (without style files

%The following (from Ams-Tex) need to be defined 
%or redefined for plain Tex

%\plainfootnote   \,   \frac   
%\text   \pmb (= Poor man's bold)   \Cal (calligraphic font)
% \comment  \endcomment   \cases \endcases
% \matrix \endmatrix

%%%% text colour definitions
\input colordvi

%%% four numbers between 0 and 1
%%% in argument of textcolour
%%% represent CMYK numbers (cyan, magenta, yellow, black

%black

%black
\comment

olive-green  7,0,95,35       green (51,0,96,39)       red (0,100,87,14)
    darkblue (94,100,0,32)   blue (98,30,0,28)      darkbrown (0,32,61,82)
    brown (0,43,83,58)      purplered (5,98,0,29)         purple (44,97,0,44)
     verydarkblue(96,97,0,84)     pink (0,28,6,4)         orange (0,53,99,31)
\endcomment
%deep blue .9 1.0 .0 .0
%deep brown .8 .8 .99 .0
% orange pastel .0 .6 .99 .0
% pink .0 .6 .3 .1
% green .9 .0 .9 .5
% dark blue .9 .7 .0 .5
%purple (sort of) .3 .7 .0 .5
% red .0 .9 .6 .3 

%\font\Ti=Times
%\def\shat{{\Ti š}}  \font\eesp=ElysiumBook at 11.3pt

%this defines eth!

%%%%%%font definitions

\font\rm=cmr10 \rm

\font\bf=cmb10
\font\Rm=cmr9 at 11pt
\rm
\font\it=cmsl9 at 10pt
 %[track+70] 
 at 7pt

\font\Rrm=cmr17 at 16pt
   \font\Rm=cmr12 at 11.5pt
%%%%%%%%%%
%\font\Frak=eufm10

  %%% (P,x^n)
\long\def\Pf{\par\noindent {\it Proof.} }
\def\({\left(}
\def\){\right)}
\def\st{such that }
\def\qed{\hfill$\bullet$\vskip 4pt}

\def\brcs#1{\left\{ #1\right\}}

\def\Log{\text{Log\,}}

\def\:{\,:}

\def\aa{\text{\bf a\,}}

\def\R{\text{\bf R}}

\def\Z{\text{\bf Z}}

\def\Arrow #1;#2.{#1\:#2 \to }
%#1: #2 -> 

\def\Set#1#2{\brcs{#1 \left|\vphantom{#1 #2} \right.#2}}
%set of #1 such that #2

%%%bigOh notation 

%%%little oh notation

%product left and right
\def\Rrr#1,#2{{\Cal J}_{#1,#2}}%difference
\def\slfrac#1#2{{\raise -.07 ex\hbox{$^{#1}$}}\!/\raise .35 ex \hbox{${}_{#2}$}}
\def\ssf #1/#2{\slfrac {#1}{#2}}

\def\pd #1,#2.{\frac {\partial #1}{\partial #2}}

%%%%%%%statement of lemmas, propositions,
%%theorems, examples, ...  \Lem 
   \long\def\Lem
#1.#2\par{\vskip4pt{\baselineskip=13pt\font\it=cmsl12 at
11.5pt\Rm
   \noindent {\rm \uppercase{#1}} #2\vskip3pt

   }} 

\long\def\Proclaim #1.#2 \endproclaim{\vskip4pt{\baselineskip=13pt\font\it=cmsl12 at
11.5pt\Rm
   \noindent {\rm \uppercase{#1}} #2\vskip3pt

   }} 

\long\def\remark #1\endremark{\vskip 2pt \noindent {\it Remark\/} #1\par}

\long\def\Sectionhead #1.#2:\par #3{\vskip 4pt \noindent {\bf #1 #2}vskip 2pt\noindent\nospace #3}

\long\def\Title #1\par {\noindent{\Rrm #1}\vskip 9pt}
%main title

 \long\def\SubT #1.{\noindent {\it #1\/} } 
%subtitle
 
 \long\def\SecT
#1\par{\vskip 3pt \noindent {\bf #1}\vglue1pt
   \noindent}%section title

\long\def\subtitle #1.{\vskip 2pt \noindent {\it #1}}

\long\def\Rmk#1\par{\vskip 1pt \noindent {\it
Remark.} #1\vskip2pt}%remark

\long\def\Abstract #1\par{{\leftskip= 3 true cm \rightskip = 3 true cm \font\it=cmsl10 \font\rm=cmr10 \baselineskip = 10pt
\parindent=.35 true cm\rm\noindent 
{\it Abstract} #1\vskip 8pt

}}

\long\def\Author #1 \par{\noindent{\it #1}}
%%%%%%%%%macros for numbering
%%of propositions, lemmas,examples, theorems

%\input (1+x)m_macros
%\input diagrams
%\input generic_macrosM

 % 

 %

%{\left\{ #1\right)}
 %

\def\Leg#1,#2.{\(\frac{#1}{{#2}}\)}
\def\Ip #1,#2.{\langle\!\langle #1,#2\rangle\! \rangle}
\def\Ipd #1,#2.{\Ip #1,#2._d}

\def\1{{\pmb 1}}

%\font\ding=Dingbats

\NoBlackBoxes

\def\cvx{\text{cvx\,}}

\def\aa{\pmb a}

\Title A Shapley-Folkman lemma for lattice polytopes

\Abstract \ Some occurrences of $n$ can be replaced by $n-1$ in a special case of the Shapley-Folkman lemma.

%%%%%%%%%%%%%%%%%
%%tensor products of dimension groups

\noindent {\it David Handelman}%
\comment
\plainfootnote{$^1$}{\hglue -.5em Supported in part by a Discovery grant from NSERC.}
\endcomment 

{}\plainfootnote{}{AMS(MOS) classification: 52B20, 52C07; key words and phrases:  lattice simplex, empty polytope, Shapley-Folkman lemma, solid polytope} 

\vskip 10pt 
\noindent The Shapley-Folkman lemma (not to be confused with the S-F {\it theorem,} which is a consequence of the lemma) asserts the following. Let $\brcs{S_i}_{i=1}^m$ be  a collection of $m$ subsets of $\R^n$, with $n < m$, and let $K_i$ denote the convex hull of $S_i$. Given $x $ in the Minkowski sum, $\sum K_i$ (the set of sums $\sum k_i$ with each $k_i \in K_i$), there exists a subset $T \subset \brcs{1,2,\dots, m}$ with $|T| = n$ \st 
$$
x = z + \sum_{s \in T^c} y_s,
$$
where $z \in \sum_T K_t$, and each $y_s \in S_s$ with $s $ running over $T^c$. Moreover, the cardinality of $T$ cannot in general be reduced to $n-1$. 

If we apply this to the special case $S =S_1 = S_2 = \dots = S_m \subset \Z^n$ with $S$ finite,  and $K = \cvx S$, then we have a result for lattice polytopes (a {\it lattice polytope\/} in $\R^n$ is the convex hull of a finite set of lattice points (elements of $\Z^n$)). Specifically, 
$$
mK \cap \Z^n = n K \cap \Z^n + (m-n)(K \cap \Z^n).  \tag 1
$$
To translate this to English (or any other language), we are using the convention that if $A$ and $B$ are subsets of an abelian group, then $A+ B$ is the set of sums $\Set{a+b}{a \in A; b \in B}$, and if $k$ is a positive integer, then $kA$ denotes $A + A + \dots + A$ ($k$ times), that is, the set of sums of $k$ elements of $A$ (despite the notation, $kA$ is not generally the same as the set of multiples by $k$ of elements of $A$, unless $A$ is convex). So what this says is that if $x$ belongs to the convex hull of $mS$ (which is $mK$) and is a lattice point, then it admits a decomposition $x = z + y$, where $y \in (m-n)S$ and $z$ is a lattice point in  $nK$. This follows immediately from the statement of the Shapley-Folkman lemma (Theorem 1 below): we obtain a decomposition $x = z + \sum_{s=1}^{m-n} y_s$ with $y_s \in S$ for all $s$, and $z \in nK$. But then $z = x - \sum_{s=1}^{m-n} y_s$, so is itself a lattice point. This argument was given by me in answer to a {\it MathOverflow\/} question, [MO].

In Theorem 3, we obtain a proof of a result that seems to have been in the \ae ther, but for which I cannot find a reference. Specifically, all occurrences of nonsuperscripted $n$ in (1) can be replaced by  $n-1$. This despite the fact that the original S-F result cannot be so improved. In other words, we will show that if $m \geq n-1$, and $S$ is a finite subset of $\Z^n$, and $K = \cvx S$, then  
$$
mK \cap \Z^n = (n-1)K \cap \Z^n + (m-n+1)(K \cap \Z^n).  \tag 2
$$
 Also, $n-1$ is sharp, in the obvious sense. We will also look at some consequences for lattice polytope classification (as in [H1, H3]) using S-F and its slight improvement.

\Lem  Theorem 1 (Shapley-Folkman Lemma [wi]). Let $\brcs{S_i }_{i=1}^m$ be  subsets of $\R^n$ with $ m> n$. For each $x \in \cvx \sum_{i=1}^m S_i$, there exists an $n$-element subset $T \subset \brcs{1,2,\dots, m}$, an element $x _0 \in \cvx \sum_T S_i$, and $x_i \in S_i$ for all $i \not\in T$ \st $x = x_0 + \sum_{T^c}  x_i$. 

\Rmk Consistent with our notation, $\sum_{i=1}^m S_i$ means the set of sums of $m$ elements, one from each of the $S_i$. 

A lattice polytope  is {\it empty\/} if it contains no lattice points other than its extreme points. If a lattice polytope is a simplex, it will be referred to as a {\it lattice simplex.} If $K$ is a lattice polytope in $\R^n$, then it is empty precisely when $K \cap \Z^n$ consists of the extreme points of $K$,  the set of which we denote $\partial_e K$. A polytope $K \subset \R^n$ {\it has interior\/} if it contains a nonempty  open ball of $\R^n$. 

\Lem Lemma 2. Let $K \subset \R^n$ be an empty lattice simplex,  and suppose that  $n \geq 2$. Then for all $x \in nK \cap \Z^n$, there exists $z \in K\cap \Z^n$ \st $x - z \in (n-1)K$. 

\Rmk Alternatively, $nK \cap \Z^n = (n-1)K \cap \Z^n + K \cap \Z^n$. 

\Pf We reduce immediately to the case that $K$ contain interior (that is, $K$ has $n+1$ extreme points) by reducing to the affine span.  By translation, we may assume the extreme points are $\brcs{v_1, v_2, \dots, v_n; 0}$. Write $x = \sum \alpha_i v_i$ with $0 \leq \alpha_i $ for all $i = 1,2,\dots, n$ and $\sum \alpha_i \leq n$. If $\sum \alpha_i \leq n-1$, then $z = 0$ will suffice. 

So we may assume that $\sum \alpha_i > n-1$. If any of the $\alpha_i$ are at least as large as $1$, then we can take $z$ to be the corresponding $v_i$. So may assume that all $\alpha_i$ are less than $1$. 

Set $y = \sum (1-\alpha_i)v_i$. As  $\sum (1-\alpha_i) = n - \sum \alpha_i < 1$, we have that $y \in K$. Since $x + y = \sum v_i \in \Z^n$, we have that $y \in \Z^n$. Hence  $y \in K \cap \Z^n = \partial_e K$. Thus $y$ is one of the extreme points, say $v_j$. But the barycentric decomposition is unique (since $K$ is a simplex), and thus $\alpha_j = 0$ and all the other $\alpha_i$ are $1$. This entails $\sum \alpha_i = n-1$, a contradiction. \qed 

\Lem Theorem 3. Let $K$ be a lattice polytope in $\R^n$, and let $m \geq n-1$ be an integer. For all $x \in mK \cap \Z^n$, there exist $x_j \in K \cap \Z^n$, $j = 1, 2, \dots, m-n+1$ and $x_0 \in (n-1)K \cap \Z^n$ \st $x = \sum_{j=0}^{m-n+1} x_j$. 

\Rmk  Alternatively, $mK \cap \Z^n = (n-1)K \cap \Z^n + (m- n+1)(K \cap \Z^n)$. 

\Pf Let $S = K \cap \Z^n$, so that $K = \cvx S$. The argument given in the second paragraph of this paper yields, via the S-F lemma, that if $x \in mK$, there is a decomposition $x = z + \sum_{1}^{m-n} y_s$ where $z \in nK \cap \Z^n$ and all $y_s \in S$. 
Hence it suffices to prove the result in the special case that $m = n$, i.e., $x \in nK \cap \Z^n$ entails there exists $x_1 \in K \cap \Z^n$ \st $x - x_1 \in (n-1)K \cap \Z^n$. 

To this end, we may triangulate $K$ by empty $n$-simplices in the obvious way; so $K = \cup K_\alpha$ where each $K_\alpha$ is an empty $n$-simplex. Then $x/n \in K$, so there exists $\alpha$ \st $x/n \in K_\alpha$, and so $x \in nK_\alpha \cap \Z^n$. By Theorem 3, we can decompose $x = x_1 + z$ where $x_1 \in K_\alpha \cap \Z^n$ and $z \in (n-1)K_\alpha \subset (n-1)K$.  As $z$ is a difference of lattice points, $z \in \Z^n$ as well, and we are done.\qed  

This is similar to, but not identical with, an old/recent result of Victoria Powers and Bruce Reznick [PR], which in connection with quadratic forms, shows that lattice points in $2K$ can be represented as a sum of distinct lattice points each in $K$, if $K = mK'$ for some integer $m \geq n-1$ and lattice polytope $K'$ ([R] uses projective dimension, so some translation will be necessary). 

It is relatively easy to construct examples showing the result is sharp, at least in one sense. In the following, $e_i$ denotes the $i$th standard column (and we only use the first $n-1$ of them). 
 
\Lem Example 4.  Let $n \geq 4$. Let $\aa = (a(i);d)\in  \Z^{1 \times n}$ with $i = 1, 2, \dots, n-1$; $0 < a(i) < d$; and $\sum a(i) <d$. Form $S = \brcs{0; e_1, e_2, \dots, e_{n-1}; \aa^T}$ and set $K = \cvx S$ in $\R^{n\times 1} = \R^n$. Then $K$ is a lattice simplex with interior, and  the column $\1 = (1,1,\dots ,1)^T \in \Z^n$ is an element of $(n-1)K \cap \Z^n$ that  cannot be represented in the form $ \sum x_i$ where $x_i \in (n-2)K \cap \Z^n$. 

\Rmk In particular, since $0 \in K $, we have $aK\cap \Z^n \subset bK \cap \Z^n$ whenever $a \leq b$. Thus there are no decompositions such as $nK \cap \Z^n = (n-i)K \cap \Z^n + i K\cap \Z^n$ if $2 \leq i \leq n-2$, or $(n-1)K\cap \Z^n = (n-2)K\cap \Z^n + K\cap \Z^n$\. This can be rephrased in terms of lattice cones. Form $\cup_{i=1}^{\infty} \(iK \cap \Z^n\)$. Then $\1$ is at atom in this cone (a nonzero minimal element), and it does not belong to any $iK$ for $i <n-1$. 

\Rmk We can easily arrange that $K$ be empty as well; for example, this holds if $1 \in \brcs{a(i)}$ and the content of the remaining $a(i)$s is $1$.

\Pf We  see that $\1 = \sum (1-a(i)/d)e_i + \aa/d$, and this decomposition (modulo the coefficient of zero) is unique, since $K$ is a simplex. The sum of the coefficients is $n-1 - (\sum a(i) - 1)/d > n-2$ (since $\sum a(i) < d$). Thus $\1 \in (n-1)K \setminus (n-2)K$. We see that if a lattice point in the cone $\cup_i \(iK \cap \Z^n\) \subset (\Z^n)^+$ has its bottom coordinate nonzero, then all of its coordinates must be strictly positive. Hence $\1$ is an atom in the cone, and the result follows from $\1$ not belonging to any $iK$ with $i < n-1$. \qed

\SecT    Local solidity  

A lattice polytope $K$ in $\R^n$ is {\it projectively faithful\/} [H1] if the set of
differences $K \cap \Z^n - K \cap \Z^n$ generates the standard copy of $\Z^n$ as
an abelian group. We say that $K$ is {\it solid\/} (op\.cit.) if for all $m$, $mK \cap \Z^n
= m(K \cap \Z^n)$ (it suffices that this be true for all $m \leq n-1$, by an application of Theorem\, 3). 

Let $v$ be an extreme point of the projectively faithful lattice polytope $K$,
and let $K_v = K - v$, that is, translate $K$ so that $v$ goes to $0$. Let $C_v$
denote the lattice cone $\cup_m m(K_v \cap \Z)$. We say  $C_v$ is {\it
unperforated\/} (or {\it normal,} or {\it integral,} \dots) if whenever  $x \in \Z^n$ and there
exists a positive integer  $m$ \st $mx \in C_v$, then $x \in C_v$ ($C_v $ has no
holes, as a subset of $\Z^n$). If for every extreme point $v$, $C_v$ has no
holes, then $K$ is {\it locally solid\/} (op\.cit.; this corresponds to the toric variety being
normal [O]). It is also equivalent to $R_P$ being integrally closed in its field of
fractions, where $P $ is a Laurent polynomial in $n$ variables with only
positive coefficients, \st $K \cap \Z^n = \Log P$ [H1; III.2 and III.8A].  In what follows, we adopt the notation of [H1].

\Lem Lemma 5. Let $K \subset \R^n$ be a projectively faithful lattice polytope.
The following are equivalent. 
\item{(a)} $K$ is locally solid; 
\item{(b)} there exists $m\geq n-1$ \st $mK \cap \Z^n = m(K \cap \Z^n)$;
 \item{(c)} for all but finitely many $m$, 
$mK \cap \Z^n = m(K \cap
\Z^n)$. 

\Pf (b) implies (c) is an immediate consequence of   Theorem 3. 

\noindent (c) implies (a). Let $v$ be a vertex, and form $C_v$; we may obviously
assume $v = 0$ (only to simplify notation); then $K \subset 2K \subset 3K
\subset \dots$. If $x $ belongs to  the convex hull of $C_v$, then $x \in mK$
for some $m$. If $m \geq n-1$, then $x \in m(K \cap \Z^n)$ by hypothesis; if $m <
n-1$, then $mK \subset (n-1)K$, and so $x \in (n-1)(K \cap \Z^n)$. So in either case, $x
\in C_v$. 

\noindent (a) implies (b). We may, as usual, assume that $0 \in K$. Form $R_P$
[H1; section I] where $P = \sum_{w \in K \cap \Z^d} x^w$. As  $K$ is locally solid, it
follows  ([H1; III.2]) that $R_P$ is integrally closed in its field of
fractions. 

Given $w \in pK \cap \Z^n$, the element $x^w/P^p$ is thus integral: by
Carath\'eodory's theorem, there exists an affinely independent subset
$\brcs{v(i)}_{i=0}^{n}$ of $\partial_e K $ \st $w = \sum \lambda_i v_i$ with
$\lambda_i \geq 0$ for all $i$, and $\sum \lambda_i = p$; from affine
independence, it follows that the solution is unique, and since all of $v_i$ and
$w$ are lattice points, all the $\lambda_i$ are rational. Hence there exists a
positive integer $a$ \st all $a\lambda_i$ are nonnegative integers (some may be
zero). Hence $(x^w/P^p)^a \in R_P$, so $x^w/P^p$ is integral over $R_P$. 

By [H1; III.2], $x^w/P^p \in R_P$. But this entails the existence of a positive
integer $N_w$ \st $w +  N_w(K\cap Z^n) \subset (p+N_w)(K\cap \Z^n)$, and  it
follows that this holds with $N_w$ replaced  by any larger integer. Let $m =p +
\max_{w \in mK \cap \Z^n} N_w$. Then for all $w \in pK \cap \Z^n$, we have  $w +
m(K\cap Z^n) \in m(K \cap \Z^n)$. Now set $p = n-1$.  

Pick $z \in m K \cap \Z^n$. We may write $z = w + (m-(n-1))(K \cap \Z^n)$ by
Theorem 3; this yields $z \in m(K \cap \Z^n)$ immediately. \qed 

\noindent{\it Two questions\/ } Let $K$ be a projectively faithful lattice
polytope.  \item{1} If $K$ is locally solid, then is it solid? 
\item{2} If $mK \cap \Z^n = m(K \cap \Z^n)$ for some $m$ \st $n/2 < m < n-1$,
then is $K$ locally solid?

\vskip 2pt

\noindent The first was suggested in [H1; p\,35]. For the second, there is a class of weird examples
given in [H3; p\,180ff]. For every $n \geq 6$, there exists a
projectively faithful lattice polytope $K \subset \R^n$ \st $iK \cap \Z^n = i(K\cap \Z^n)$
for all $i \leq n/2$, but for no $i > n/2$.

\long\def\Rf[#1] #2, #3. #4\par%
{\vskip 2pt \itemitem{[#1]} #2, {\it #3,} #4\par\vskip2pt}

\SecT References

\Rf [H1] D Handelman, Positive polynomials, convex polytopes, and a random walk problem. Lecture Notes in Mathematics 1082 (1986) Springer-Verlag 133 + x.

\Rf [H2] D Handelman, Positive polynomials and product type actions of compact groups. Memoirs of the American Mathematical Society  310 (1985)  79 + xi.

\Rf [H3] D Handelman, Effectiveness of an affine invariant for indecomposable integral polytopes. J Pure \& App Algebra  66 (1990)  165--184.

\Rf [MO] D Handelman, https:/\!/mathoverflow.net/questions/204200/integer-decomposition-of-dilated-integral-polytopes/321660\#321660.  MathOverFlow answer.

\Rf [O] T Oda, Convex bodies and algebraic geometry{\/\rm:} An introduction to toric varieties.  Ergebnisse der Mathematik und ihrer Grenzgebiete 3, Folge, Band 15, Springer-Verlag, Berlin, Heidelberg, New York, 1988.

\Rf [PR] V Powers \& B Reznick, A note on mediated simplices. https://arxiv.org/pdf/1909.11008.pdf

\Rf [Wi] Shapley-Folkman lemma, https:/\!/en.wikipedia.org/wiki/Shapley\%E2\%80\%93Folkman$\_$lemma. Wikipedia.

{}

\vskip 10pt

Mathematics Department, University of Ottawa, Ottawa ON  K1N 6N5, Canada; dehsg\@uottawa.ca

 \end